\documentclass[12pt]{article}
\usepackage{mathrsfs}
\usepackage{epic,eepic,epsf,epsfig}
\usepackage{amsfonts,srcltx,mathrsfs}
\usepackage{multirow}
\usepackage{amsmath}
\usepackage{amssymb}
\usepackage{amsbsy}
\usepackage{graphicx}
\usepackage{amsfonts}
\usepackage{color}
\usepackage{setspace}

\newtheorem{lem}{Lemma}[section]%
\newtheorem{theorem}[lem]{Theorem}%
\newtheorem{exam}[lem]{Example}%
\newtheorem{prop}[lem]{Proposition}%

\def\a{\alpha}

 \def\O{\Omega} \def\G{\Gamma}

\def\di{\bigm|}  
\def\nd{\mathrel{\bigm|\kern-.7em/}}

\def\f{\noindent}

\def\PSL{\hbox{\rm PSL}}\def\PSU{\hbox{\rm PSU}}
  \def\GL{\hbox{\rm GL}} \def\Mult{\hbox{\rm Mult}}
\def\PSp{\hbox{\rm PSp}}\def\P\GammaL{\hbox{\rm P\Gamma L}} 
 
\def\AGL{\hbox{\rm AGL}}

\def\Aut{\hbox{\rm Aut}}

\def\Cay{\hbox{\rm Cay}}
\def\Cos{\hbox{\rm Cos}}

\newcommand{\qed}{\mbox{\raisebox{0.7ex}{\fbox{}}} \vspace{4truemm}}
\def\mz{{\mathbb Z}}

\begin{document}
\title{Tetravalent $2$-arc-transitive Cayley graphs on non-abelian simple groups}

\author{Jia-Li Du, Yan-Quan Feng\footnotemark\\
{\small\em Department of Mathematics, Beijing
Jiaotong University, Beijing 100044, China}}

\footnotetext[1]{ Corresponding author. E-mails:
JiaLiDu@bjtu.edu.cn, yqfeng$@$bjtu.edu.cn}

\date{}
 \maketitle

\begin{abstract}

A graph $\G$ is said to be {\em $2$-arc-transitive} if  its full automorphism group $\Aut(\G)$ has a single orbit on ordered paths of length $2$, and for $G\leq \Aut(\G)$, $\G$ is {\em $G$-regular} if $G$ is regular on the vertex set of $\G$. Let $G$ be a finite non-abelian simple group and let $\G$ be a
connected tetravalent $2$-arc-transitive $G$-regular graph. In 2004, Fang, Li and Xu proved that either $G\unlhd \Aut(\G)$ or $G$ is one of $22$ possible candidates. In this paper, the number of candidates is reduced to $7$, and for each candidate $G$, it is shown that $\Aut(\G)$ has a normal arc-transitive non-abelian simple subgroup $T$ such that $G\leq T$ and the pair $(G,T)$ is explicitly given.

\bigskip
\f {\bf Keywords:} Cayley graph, coset graph, simple group.

\end{abstract}

\section{Introduction}

Throughout this paper, all groups and graphs are finite, and all graphs are simple and undirected. Let $G$ be a permutation group on a set $\O$ and let $\a\in \O$. Denote by $G_{\a}$ {\it the stabilizer} of $\a$ in $G$, that is, the subgroup of $G$ fixing the point $\a$. The group $G$ is {\it semiregular} if $G_{\a}=1$ for every $\a\in\O$, and {\it regular} if $G$ is transitive and semiregular.

For a graph $\G$, denote by $V(\G)$, $E(\G)$ and $\Aut(\G)$ its vertex set, edge set and full automorphism group, respectively. An {\it $s$-arc} in $\G$ is an ordered $(s+1)$-tuple $(v_0,v_1,...,v_s)$ of vertices of $\G$ such that $v_{i-1}$ is adjacent to $v_i$ for $1\leq i\leq s$, and $v_{i-1}\neq v_{i+1}$ for $1\leq i<s$.  Let $G\leq \Aut(\Gamma)$. The graph $\G$ is said to be {\it $(G,s)$-arc-transitive} or  {\it $G$-regular} if $G$ has a single orbit on $s$-arcs of $\G$ or is regular on vertices of $\G$, respectively. For short, a $1$-arc means an {\it arc}, and  $(G,1)$-arc-transitive means {\it $G$-arc-transitive}. The graph $\G$ is said to be {\it $s$-arc-transitive} if it is $(\Aut(\G), s)$-arc-transitive. In particular, $0$-arc-transitive is {\it vertex-transitive}, and $1$-arc-transitive is {\it arc-transitive} or {\it symmetric}.

For a group $G$ and a subset S of $G$ such that $1\notin S$ and $S^{-1}=S$, the \emph{Cayley graph} $\Cay(G,S)$ on $G$ with respect to $S$ is defined to have vertex set $G$ and edge set $\{\{g, sg\} | g \in G, s \in S\}$. For $g\in G$, the map $R(g): x\mapsto xg$ for $x\in G$ is a permutation on $G$, and $R(G)=\{R(g)\ |\ g\in G\}$ consists of a permutation group on $G$, called the right regular representation of $G$. It is easy to see that $R(G)\leq \Aut(\G)$. A Cayley graph $\G=\Cay(G,S)$ is said to be \emph{normal} if $R(G)$ is a normal subgroup of $\Aut(\G)$, and in this case, $\Aut(\G)=R(G)\rtimes \Aut(G,S)$ by Godsil~\cite{Godsil} or Xu~\cite{Xu}, where $\Aut(G,S)=\{\alpha \in \Aut(G)\ |\ S^{\alpha}=S\}$. Note that $\Cay(G,S)$ is $R(G)$-regular. It is well-known that if $\G$ is $G$-regular then $\G$ is isomorphic to a Cayley graph on $G$.

The investigation of symmetric graphs has a long and interesting history, highlighted at an early stage by ingenious work by Tutte~\cite{Tutte1,Tutte2} on the cubic case. Let $\G$ be a connected symmetric cubic Cayley graph on a non-abelian simple group $G$. Li~\cite{CHLi} proved that either $\G$ is normal or $G=A_5$, $A_7$, $\PSL(2,11)$, $M_{11}$, $A_{11}$, $A_{15}$, $M_{23}$, $A_{23}$ or $A_{47}$. Based on Li's work, Xu {\em et al}~\cite{XFWX2005,XFWX} proved that either $\G$ is normal or $G= A_{47}$, and for the latter, $\G$ is $5$-arc-transitive and there are only two such graphs.

Let $\G$ be a connected tetravalent $2$-arc-transitive Cayley graph on a non-abelian simple group $G$. Fang {\em et al}~\cite{Fang4} proved that either $\G$ is normal, or $G$ is one of the $22$ possible candidates listed in~\cite[Table~1]{Fang4} (the group $A_{23}$ is missed in the table). The following is the main result of this paper.

\begin{theorem}\label{theo=main}
Let $G$ be a non-abelian simple group and $\G$ a connected tetravalent $2$-arc-transitive $G$-regular graph. Then either $G\unlhd \Aut(\G)$ or $\Aut(\G)$ contains a normal arc-transitive non-abelian simple subgroup $T$ such that $G\leq T$ and $(G,T)$ is listed in Table~\ref{table=1}.
\begin{table}[ht]
\begin{center}
\begin{tabular}{|c|c|c|c|c|c|c|c|}

\hline
$G$  & $M_{11}$ & $A_{2^3\cdot 3-1}$  &$A_{2^2\cdot 3^2-1}$  & $A_{2^3\cdot 3^2-1}$  & $A_{2^4\cdot 3^2-1}$  &$A_{2^4\cdot 3^3-1}$  &$A_{2^4\cdot 3^6-1}$    \\
\hline
$T$  & $M_{12}$ & $A_{2^3\cdot 3}$   & $A_{2^2\cdot 3^2}$ & $A_{2^3\cdot 3^2}$  & $A_{2^4\cdot 3^2}$ & $A_{2^4\cdot 3^3}$  & $A_{2^4\cdot 3^6}$   \\
\hline
\end{tabular}
\end{center}
\vskip -0.5cm
\caption{{7 possible pairs of non-abelian simple groups}}\label{table=1}
\end{table}
\end{theorem}

For connected  symmetric cubic Cayley graphs on non-abelian simple groups, similar to Theorem~\ref{theo=main}, there are six possible pairs $(G,T)=(A_{47},A_{48})$, $(\PSL(2,11),M_{11})$, $(M_{11},M_{12})$, $(A_{11}, A_{12})$, $(M_{23},M_{24})$ or $(A_{23},A_{24})$ (see \cite[Theorem 7.1.3]{CHLi}), and Xu {\em et al}~\cite{XFWX2005, XFWX} proved that only the pair $(G,T)=(A_{47},A_{48})$ can happen and there are exactly two connected non-normal symmetric cubic Cayley graphs on $A_{47}$. For the $7$ possible pairs of $(G,T)$ in Theorem~\ref{theo=main}, by MAGMA~\cite{magma} there is only one $2$-arc-transitive Cayley graph on $M_{11}$ for $(G,T)=(M_{11},M_{12})$ (see  Remark of Lemma~\ref{lem=insolvable}), and there are four $2$-arc-transitive Cayley graphs on $A_{23}$ for $(G,T)=(A_{23},A_{24})$. The number of $3$-arc-transitive Cayley graphs on $A_{35}$ for $(G,T)=(A_{35},A_{36})$ is $4$, on $A_{71}$ for $(G,T)=(A_{71},A_{72})$ is $18$, and on $A_{143}$ for $(G,T)=(A_{143},A_{144})$ is $31$. At the end of this paper, we give examples of connected non-normal $2$-arc-transitive Cayley graph on $A_{23}$ for $(G,T)=(A_{23},A_{24})$ (see Example~\ref{NonNormalExample}).

\section{Preliminaries}

In this section, we describe some preliminary results which will be used later. The first result is the stabilizers of tetravalent $2$-arc-transitive graphs, given in \cite[Theorem 4]{Potocnik}.

\begin{prop}\label{prop=stabilizer}
Let $\Gamma$ be a connected tetravalent $(G,s)$-arc-transitive but not   $(G,s+1)$-arc-transitive graph with $v\in V(\Gamma)$. If $s\geq 2$ then one of the following occurs:
\begin{enumerate}

\itemsep -1pt
\item [\rm (1)] For $s=2$, we have $G_v\cong A_4$ or $S_4$. In particular, $|G_v|=2^2\cdot 3$ or $2^3\cdot 3$.
\item [\rm (2)] For $s=3$, we have $G_v\cong \mz_3 \times A_4$, $\mz_3\rtimes S_4$, or $S_3 \times S_4$. In particular, $|G_v|=2^2\cdot 3^2$, $2^3\cdot 3^2$ or $2^4\cdot 3^2$.
\item [\rm (3)] For $s=4$, we have $G_v\cong\mz_2^3\rtimes \rm GL(2,3)= \rm AGL(2,3)$. In particular, $|G_v|=2^4\cdot 3^3$.
\item [\rm (4)] For $s=7$, we have $G_v\cong \mz_3^5 \rtimes \GL(2,3)$. In particular, $|G_v|=2^4\cdot 3^6$.
\end{enumerate}
\end{prop}

\f {\bf Remark:} For $s=7$, by \cite[Theorem 1.1]{Potocnik}, $G_v=\langle e_0,e_1,e_2,e_3,e_4,e_5,e_6,d \rangle\cong$
$\mz_3^5\rtimes\rm GL(2,3)$ with the following relations: $e_i^3=d^2=1$ for all $0\leq i\leq 6$, $[e_i,e_j]=1$ for
all $|i-j|<4$, $[e_0,e_4]=e_2$, $[e_0,e_5]=e_1^{-1}e_2e_3^{-1}e_4^{-1}$, $e_0e_6e_0=e_6e_0e_6$, $[e_1,e_5]=e_3$,
$e_1^{e_6}=e_5e_4e_3^{-1}e_2e_1$, $e_2^{e_6}=e_2e_4^{-1}$, $e_k^d=e_k^{-1}$
for $k=0,2,3,5,6$ and $e_l^d=e_l$ for $l=1,4$, and by MAGMA~\cite{magma}, $G_v$ has no normal subgroup isomorphic to $A_4$, $S_4$, $\mz_3\times A_4$, $\mz_3\rtimes S_4$, $S_3\times S_4$, $\rm AGL(2,3)$, or a non-trivial $2$-group.

\medskip

Let $\G$ be a graph and let $N \leq \Aut(\G)$. The \emph{quotient graph} $\G_N$
of $\G$ relative to $N$ is defined as the graph with vertices
the orbits of $N$ on $V(\G)$ and with two orbits adjacent if there is an edge
in $\G$ between those two orbits. The theory of quotient graph is widely used
to investigate symmetric graphs. The following proposition can be deduced from  \cite[Theorem 1.1]{Gardiner} and \cite[Theorem 4.1]{Praeger}.

\begin{prop}\label{prop=atlesst3orbits}
Let $\Gamma$ be a connected tetravalent $(G,s)$-arc-transitive graph for some $s\geq 2$ and let $N$ be a normal subgroup of $G$. If $N$ is transitive then either $N$ is regular or arc-transitive, and if $N$ has at least three orbits then $N$ acts semiregularly on $V(\G)$ and the quotient graph $\G_N$ is a connected tetravalent $(G/N,s)$-arc-transitive graph.
\end{prop}

Guralnick~\cite{Guralnick} classified non-abelian simple groups which contain subgroups of index a power of a prime, and by \cite[Theorem 1]{LiXu}, we have the following result.

\begin{prop}\label{prop=simplegroup}
Let $G$ and $T$ be non-abelian simple groups such that $G\leq T$ and $|T:G|=2^a\cdot 3^b\geq 6$ with $0\leq a\leq 4$ and $0\leq b\leq 6$.
Then $T$, $G$ and $|T:G|$ are listed in Table~\ref{table=2}.
\end{prop}

\begin{table}[ht]
\begin{center}
\begin{tabular}{|c|c|c||c|c|c|}

\hline
$T$              & $G$              & $|T:G|$              &$T$              & $G$               & $|T:G|$         \\
\hline
$M_{11}$         & $\PSL(2,11)$     & $2^2\cdot3$          &$M_{12}$         & $M_{11}$          & $2^2\cdot 3$         \\
\hline
$M_{24}$         & $M_{23}$         & $2^3\cdot 3$         &$\PSU(3,3)$      & $\PSL(2,7)$       & $2^2\cdot3^2$         \\
\hline
$A_9$            & $A_7$            & $2^3\cdot 3^2$       &$\PSp(4,3)$      & $A_6$             & $2^3\cdot 3^2$         \\
\hline
$\PSp(6,2)$     & $A_8$            & $2^3\cdot 3^2$       &$\PSU(4,3)$      & $\PSL(3,4)$       & $2\cdot3^4$          \\
\hline
$M_{12}$        & $\PSL(2,11)$     & $2^4\cdot3^2$         & $\PSU(4,3)$      & $A_7$            & $2^4\cdot3^4$         \\
\hline
$G_2(3)$        & $\PSL(2,13)$    & $2^{4}\cdot3^5$       &$A_n$      & $A_{n-1}$                & $n=2^a\cdot3^b$      \\
\hline
\end{tabular}
\end{center}
\vskip -0.5cm
\caption{{Non-abelian simple group pairs of index $2^a\cdot 3^b$}}\label{table=2}
\end{table}

Let $G$ and $E$ be two groups. We call an extension $E$ of $G$ by $N$ a {\em central extension} of $G$ if $E$ has a central subgroup $N$ such that $E/N\cong G$, and if further $E$ is perfect, that is, the derived group $E'=E$, we call $E$ a {\em covering group} of $G$. A covering group $E$ of $G$ is called a {\em double cover} if $|E|=2|G|$. Schur~\cite{Schur} proved that for every non-abelian simple group $G$ there is a unique maximal covering group $M$ such that every covering group of $G$ is a factor group of $M$ (see \cite[Kapitel V, \S23]{Huppert}). This group $M$ is called the {\em full covering group} of $G$, and the center of $M$ is the {\em Schur multiplier} of $G$, denoted by $\Mult(G)$.

\begin{prop}\label{prop=covering group} For $n\geq 5$, the alternating group $A_n$ has a unique double cover $2.A_n$, and for $n\geq 7$,   all subgroups of index $n$ of $2.A_n$ are conjugate and isomorphic to $2.A_{n-1}$.
\end{prop}

\f {\bf Proof:} By Kleidman and Liebeck~\cite[Theorem 5.1.4]{Kleidman}, $\Mult(A_n)\cong \mz_2$ for $n\geq 5$ with $n\not=6,7$,  and $\Mult(A_n)\cong \mz_6$ for $n=6$ or $7$. This implies that $A_n$ has a unique double cover for $n\geq 5$, and we denote it by $2.A_n$. Since $A_n$ has no proper subgroup of index less than $n$, all subgroups of index $n$ of $2.A_n$ contain the center of order $2$ of  $2.A_n$. Let $n\geq 7$. By~\cite[2.7.2]{Wilson}, $2.A_n$ contains a subgroup $2.A_{n-1}$ of index $n$, and since all subgroups of index $n$ of $A_n$ are conjugate, all subgroups of index $n$ of $2.A_n$ are conjugate and hence isomorphic to $2.A_{n-1}$.
\hfill\qed

Now, we introduce the so called coset graph. Let $G$ be a group and $H\leq G$. Denote by $D$ a union of some double cosets of $H$ in $G$ such that $D^{-1}=D$. The \emph{coset graph} $ \Gamma=\Cos(G,H,D)$ on $G$ with respect to $H$ and $D$ is defined to have vertex set $V(\Gamma)=[G:H]$, the set of right cosets of $H$ in $G$, and the edge set $E(\Gamma)=\{\{Hg,Hdg\}|~g\in G,d\in D\}$. It is well known that
$ \Gamma=\Cos(G,H,D)$ has valency $|D|/|H|$ and it is connected if and only if $G=\langle D,H\rangle$, that is, $D$ and $H$ generate $G$.
The action of $G$ on $[G:H]$ by right multiplication induces a vertex-transitive group of automorphisms of $\G$, and this group is arc-transitive if and only if $D$ is a single double coset. Moreover, this action is faithful if and only if $H_G=1$, where $H_{G}$ is the largest normal subgroup of $G$ contained in $H$. Clearly, $\Cos(G,H,D)\cong \Cos(G,H^{\alpha},D^{\alpha})$ for each $\alpha \in \Aut(G)$.

Conversely, let $\Gamma$ be a $G$-vertex-transitive graph.
By \cite{Sabidussi}, $\Gamma$ is isomorphic to a coset graph $\Cos(G,H,D)$,
where $H=G_{v}$ is the vertex stabilizer of $v\in V(\Gamma)$ in $G$ and $D$
consists of all elements of $G$ which map $v$ to one of its neighbors.
It is easy to show that $H_G=1$ and $D$ is a union of some double cosets
of $H$ in $G$ satisfying $D^{-1}=D$. Assume that $G$ is arc-transitive
and $g \in G$ interchanges $v$ and one of its neighbors. Then $g^2 \in H$
and $D=HgH$. Furthermore, $g$ can be chosen as a $2$-element in $G$,
and the valency of $\Gamma$ is $|D|/|H|=|H:H\cap H^g|$. For more
details regarding coset graph, refer to \cite{FangP,Lorimer,Miller,Sabidussi}.

\begin{prop}\label{prop=cosetgraph}
Let $\Gamma$ be a connected $G$-arc-transitive graph and  $\{u,v\}$ an edge of $\Gamma$.
Then $ \Gamma $ is isomorphic to a coset graph $\Cos(G,G_{v},G_{v}gG_{v})$,
where $g$ is a $2$-element in $G$ such that $G_{uv}^{g}=G_{uv}$, $g^2\in G_{v}$
and $\langle G_v,g\rangle=G$. Moreover, $\Gamma$ has valency $|G_v:G_v \cap G_{v}^{g}|$.
\end{prop}

In Proposition \ref{prop=cosetgraph}, the $2$-element $g$ is called {\em feasible} to $G$ and $G_v$. Feasible $g$ can be found by MAGMA~\cite{magma} when the order $|G|$ is not too large, and for convenience, an example of computer program to find feasible $g$ is given as an appendix at the end of the paper, which will be used in Section~\ref{s3} frequently.

\section{Proof of Theorem~\ref{theo=main}\label{s3} }

In this section, we always assume that $G$ is a non-abelian simple group.
To prove Theorem~\ref{theo=main}, we need the following three lemmas.

\begin{lem}\label{A12}
There is no connected tetravalent $A_{12}$-arc-transitive $A_{11}$-regular graph, and no connected tetravalent $M_{12}$-arc-transitive graph with stabilizer isomorphic to $S_4$.
\end{lem}

\f {\bf Proof:} Let $\G$ be a connected tetravalent $T$-arc-transitive graph with $v\in V(\G)$. By Proposition~\ref{prop=cosetgraph},  $\G=\Cos(T,T_v,T_vtT_v)$ for some feasible $t$, that is, $t$ is a $2$-element such that $t^2\in T_v$, $\langle T_v,t\rangle=T$ and $|T_v:T_v\cap T_v^t|=4$.

Let $T=A_{12}$ and let $\G$ be $A_{11}$-regular. Then $|A_{11}|=|V(\G)|$ and $|T_v|=|T|/|A_{11}|=12$. By Proposition~\ref{prop=stabilizer}, $T_v\cong A_4$. By MAGMA~\cite{magma}, $T_v$ has $12$ conjugacy classes in $T$. Take a given $T_v$ in each conjugacy class, and computation shows that there is no such feasible $t$. It follows that there is no connected tetravalent $A_{12}$-arc-transitive $A_{11}$-regular graph.

Let $T=M_{12}$ and $T_v\cong S_4$. By MAGMA, $T_v$ has four conjugacy classes in $T$. Take a given $T_v$ in each conjugacy class, and computation shows that there is no such feasible $t$. Thus, there is no connected tetravalent $M_{12}$-arc-transitive graph with stabilizer isomorphic to $S_4$.
\hfill\qed

The {\em radical} of a group is defined as its largest soluble normal subgroup.

\begin{lem}\label{lem=insolvable}
Let $\G$ be a connected tetravalent $G$-vertex-transitive graph and let $X$ be a $2$-arc-transitive subgroup of $\Aut(\G)$ containing $G$.
If $X$ has trivial radical, then either $G\unlhd X$, or $X$ has a normal arc-transitive non-abelian simple  subgroup $T$ such that $G\leq T$ and $(G,T)=(M_{11},M_{12})$ or $(A_{n-1},A_n)$ with $n\geq 8$ and $n\di 2^4 \cdot 3^6$.
Moreover, if $\G$ is $G$-regular then $(G,T)$ is listed in Table~\ref{table=1}.
\end{lem}

\f {\bf Proof:} Let $N$ be a minimal normal subgroup of $X$. Since $X$ has trivial radical, $N=T^s$ for a positive integer $s$ and a non-abelian simple group $T$. Clearly, $NG\leq X$. Since $\Gamma$ is $G$-vertex-transitive,  by the Frattini argument we have $X=GX_v$ for $v\in V(\G)$, and hence $|X|=|G||X_v|/|G_v|$. Since $|NG|=|N||G|/|N \cap G|$, we have $|N|/|N\cap G|\di |X_v|/|G_v|$, and by Proposition~\ref{prop=stabilizer}, $|N|/|N\cap G|\di 2^4\cdot 3^6$.

The simplicity of $G$ implies that $N \cap G =1$ or $G$. If $N \cap G =1$ then $|N|=|N|/|N\cap G|\di 2^4\cdot 3^6$, which is impossible because  $N$ is insoluble. Thus, $N \cap G =G$, that is, $G\leq N$, and $|N|/|G|\di 2^4\cdot 3^6$. Since $T\unlhd N$, we have $T\cap G=1$ or $G$. If $T\cap G=1$ then
$|T|=|T|/|G\cap T|=|GT|/|G|\di |N|/|G|$, which also is impossible because $|N|/|G|\di 2^4\cdot 3^6$. Thus, $G\cap T=G$, that is, $G\leq T$. Since $|N|/|G|=|T|^{s-1}|T:G|\di 2^4\cdot 3^6$, we have $s=1$. It follows that $G\leq T\unlhd X$ and $|T:G|=|T|/|G|\di 2^4\cdot 3^6$.

If $G=T$ then $G\unlhd X$ and we are done. Now we may assume  that $G$ is a proper subgroup of $T$. In particular, $G\ntrianglelefteq X$.
By Proposition~\ref{prop=simplegroup}, $(G,T)$ is listed in Table~\ref{table=2}. Since $G\leq T\lhd X$, $\G$ is $T$-vertex-transitive, and Proposition \ref{prop=atlesst3orbits} implies that $\G$ is $T$-arc-transitive. By Proposition~\ref{prop=cosetgraph}, $\G=\Cos(T,T_v,T_vtT_v)$ for some feasible $t$.
Note that $|T|/|G|=|V(\G)||T_v|/(|V(\G)||G_v|)=|T_v|/|G_v|$.

Suppose $(G,T)=(\PSL(2,13),G_2(3))$. Then $|T_v|/|G_v|=|T|/|G|=2^4\cdot3^5\di |T_v|$. By Proposition~\ref{prop=stabilizer}, $T_v\cong \mz_3^5\rtimes \GL(2,3)$, and $|V(\G)|=|T|/|T_v|=364$. However, there is no connected tetravalent $2$-arc-transitive graph of order $364$ by \cite[Table 2]{Potocnik}, a contradiction.

Suppose $(G,T)=(\PSL(3,4),\PSU(4,3))$ or $(A_7,\PSU(4,3))$. Then  $3^4\di |T_v|$ because $|T_v|/|G_v|=|T|/|G|$. By Proposition~\ref{prop=stabilizer}, $T_v\cong \mz_3^5\rtimes \GL(2,3)$, and by MAGMA, $\PSU(4,3)$ has no subgroup isomorphic to $\mz_3^5\rtimes \GL(2,3)$, a contradiction.

Suppose $(G,T)=(A_6,\PSp(4,3))$. Then $|T_v|/|G_v|=|T|/|G|=2^3 \cdot 3^2\di |T_v|$. Since $|T|=2^6 \cdot 3^4 \cdot5$, Proposition~\ref{prop=stabilizer} implies that $T_v\cong \mz_3\rtimes S_4$, $S_3 \times S_4$ or $\rm AGL(2,3)$, and by MAGMA, $\PSp(4,3)$ has no subgroup isomorphic to $\mz_3\rtimes S_4$, $S_3 \times S_4$ or $\rm AGL(2,3)$, a contradiction.

Suppose $(G,T)=(\PSL(2,7),\PSU(3,3))$. Then $|T_v|/|G_v|=|T|/|G|=2^2 \cdot 3^2\di |T_v|$. Since $|T|=2^5 \cdot 3^2 \cdot7$, Proposition~\ref{prop=stabilizer} implies that $T_v\cong \mz_3\times A_4$, $\mz_3\rtimes S_4$ or $S_3 \times S_4$, and by MAGMA, $\PSU(3,3)$ has no subgroup isomorphic to $\mz_3\times A_4$, $\mz_3\rtimes S_4$ or $S_3 \times S_4$, a contradiction.

Suppose $(G,T)=(\PSL(2,11),M_{11})$. Then $|T_v|/|G_v|=|T|/|G|= 2^2\cdot 3\di |T_v|$. Since $|T|=2^4\cdot 3^2\cdot 5\cdot 7$,  Proposition~\ref{prop=stabilizer} implies that $T_v\cong A_4 $, $S_4$, $\mz_3\times A_4$, $\mz_3\rtimes S_4$ or $S_3\times S_4$. By MAGMA, $M_{11}$ has no subgroup isomorphic to $\mz_3\times A_4$, $\mz_3\rtimes S_4$ or $S_3\times S_4$, and if $T_v\cong A_4 $ or $S_4$ then $T_v$ has one conjugacy class, respectively. By taking a given $T_v$ in each conjugacy class, computation shows that there is no feasible $t$, a contradiction.

Suppose $(G,T)=(A_7,A_9)$. Then $|T_v|/|G_v|=|T|/|G|=2^3\cdot 3^2\di |T_v|$. Since $3^5\nmid |A_9|$,  Proposition~\ref{prop=stabilizer} implies that
$T_v\cong \mz_3\rtimes S_4 $, $S_3\times S_4$ or $\rm AGL(2,3)$. By MAGMA, $A_9$ has no subgroup isomorphic to $\rm AGL(2,3)$ and if $T_v\cong \mz_3\rtimes S_4 $ or $S_3\times S_4$ then $T_v$ has two conjugacy classes, respectively. By taking a given $T_v$ in each conjugacy class, computation shows that there is no feasible $t$, a contradiction.

Suppose $(G,T)=(M_{23},M_{24})$. Then $|T_v|/|G_v|=|T|/|G|= 2^3\cdot3$ and $2^3\cdot3\di |T_v|$. Since $3^4\nmid |M_{24}|$,  Proposition~\ref{prop=stabilizer} implies that $T_v\cong S_4$, $ \mz_3\rtimes S_4 $, $S_3\times S_4$ or $\rm AGL(2,3)$. By MAGMA, if $T_v\cong S_4$ then $T_v$ has $19$ conjugacy classes, if $T_v\cong \rm AGL(2,3)$ then $T_v$ has one conjugacy class, and if $T_v\cong \mz_3\rtimes S_4 $ or $S_3\times S_4$ then $T_v$ has four conjugacy classes, respectively. By taking a given $T_v$ in each conjugacy class, computation shows that there is no feasible $t$, a contradiction.

Suppose $(G,T)=(A_8,\PSp(6,2))$. Then $|T_v|/|G_v|=|T|/|G|= 2^3\cdot3^2\di |T_v|$. Since $3^5\nmid |\PSp(6,2)|$,  Proposition~\ref{prop=stabilizer} implies $T_v\cong \mz_3\rtimes S_4 $, $S_3\times S_4$ or $\rm AGL(2,3)$. By MAGMA, $\PSp(6,2)$ has no subgroup isomorphic to $\rm AGL(2,3)$. If $T_v\cong \mz_3\rtimes S_4 $ then $T_v$ has four conjugacy classes, and if $T_v\cong S_3\times S_4$ then $T_v$ has eight conjugacy classes. By taking a given $T_v$ in each conjugacy class, computation shows that there is no feasible $t$, a contradiction.

Suppose $(G,T)=(\PSL(2,11),M_{12})$. Then $|T_v|/|G_v|=|T|/|G|= 2^4\cdot3^2\di |T_v|$. Since $|T|=2^6 \cdot 3^3 \cdot5\cdot 11$, Proposition~\ref{prop=stabilizer} implies that $T_v\cong S_3\times S_4$ or $\rm AGL(2,3)$. By MAGMA, $M_{12}$ has no subgroup isomorphic to $S_3\times S_4$ and if $T_v\cong \rm AGL(2,3) $ then $T_v$ has two conjugacy classes. By taking a given $T_v$ in each conjugacy class, computation shows that there is no feasible $t$, a contradiction.

Suppose $(G,T)=(A_5,A_6)$ Then  $|T_v|/|G_v|=|T|/|G|=2\cdot 3\di |T_v|$. By Atlas~\cite[pp.4]{Atlas} and Proposition~\ref{prop=stabilizer}, $T_v\cong A_4$ or $S_4$, and by MAGMA, $T_v$ has two conjugacy classes, respectively. By taking a given $T_v$ in each conjugacy class, computation shows that there is no feasible $t$, a contradiction.

By the above contradictions, $(G,T)\not=(\PSL(2,13),G_2(3))$, $(\PSL(3,4)$, $\PSU(4,3))$, $(A_7$, $\PSU(4,3))$, $(A_6$, $(\PSp(4,3))$, $(\PSL(2,7)$, $\PSU(3,3))$, $(\PSL(2,11)$, $M_{11})$, $(A_7$, $A_9)$, $(M_{23}$, $M_{24})$, $(A_8$, $\PSp(6,2))$, $(\PSL(2,11)$, $M_{12})$ or $(A_5$, $A_6)$. Deleting the above impossible pairs from Table~\ref{table=2}, we have $(G,T)=(M_{11},M_{12})$ or $(A_{n-1},A_n)$ with $n\geq 8$ and $n\di 2^4\cdot 3^6$ because $|T|/|G|\di 2^4\cdot 3^6$. To finish the proof, let $\G$ be $G$-regular. Then $G_v=1$ and  $|T_v|=|T_v|/|G_v|=|T|/|G|$. By Proposition \ref{prop=stabilizer}, $(G,T)=(M_{11},M_{12})$ or $(A_{n-1},A_n)$ with $n=2^2\cdot 3$, $2^3\cdot 3$, $2^2\cdot3^2$, $2^3\cdot3^2$, $2^4\cdot3^2$, $2^4\cdot3^3$ or $2^4\cdot3^6$, and by Lemma~\ref{A12}, $(G,T)\not=(A_{11},A_{12})$. It follows that $(G,T)$ is listed in Table~\ref{table=1}.
\hfill\qed

\f {\bf Remark:} Let $(G,T)=(M_{11},M_{12})$. Then there is a unique connected tetravalent $T$-arc-transitive $G$-regular graph $\G$, and $\Aut(\G)\cong M_{12}\rtimes\mz_6$ has  non-trivial radical $\mz_3$. These facts can be checked by MAGMA. In fact, since $|T_v|=|T_v|/|G_v|=|T|/|G|= 2^2\cdot3$, Proposition~\ref{prop=stabilizer} implies that $T_v\cong A_4$.  By Proposition~\ref{prop=cosetgraph}, $\G=\Cos(T,T_v,T_vtT_v)$ for some feasible $t$, and by MAGMA, computation shows that $T_v$  has four conjugacy classes in $T$. Take a given $T_v$ in each conjugacy class: for two conjugacy classes there is no feasible $t$, and for the other two conjugacy classes, one has $24$ feasible $t$ but all corresponding graphs $\G$ are not $M_{11}$-vertex-transitive, and the other has $12$ feasible $t$, of which the corresponding graphs $\G$ are isomorphic to each other and $\Aut(\G)=M_{12}\rtimes \mz_6$ with radical $\mz_3$.
\hfill\qed

\begin{lem}\label{lem=T normal}
Let $\G$ be a connected tetravalent $2$-arc-transitive $G$-regular graph and let $\Aut(\G)$ have non-trivial radical $R$ with at least three orbits on $V(\G)$. Assume $RG=R\times G$. Then $G\unlhd \Aut(\G)$ or $\Aut(\G)$ contains a normal arc-transitive non-abelian simple  subgroup $T$ such that $G\leq T$ and $(G,T)$ is listed in Table~\ref{table=1}.
\end{lem}

\f {\bf Proof:} Set $A=\Aut(\G)$ and $B=RG=R\times G$. Then $G$ is characteristic in $B$. To finish the proof, we may assume $G\ntrianglelefteq A$ and aim to show that $A$ contains a normal arc-transitive non-abelian simple  subgroup $T$ such that $G\leq T$ and $(G,T)$ is listed in Table~\ref{table=1}.

 Since $R\neq 1$ has at least three orbits, by Proposition~\ref{prop=atlesst3orbits} the quotient graph $\G_{R}$ is a connected tetravalent $(A/R,2)$-arc-transitive graph with $A/R\leq \Aut(\G_{R})$, and $\G_R$ is $B/R$-vertex-transitive. Since $G\ntrianglelefteq A$ and $G$ is characteristic in $B$, we have $B\ntrianglelefteq A$ and $G\cong B/R\ntrianglelefteq A/R$. Furthermore, $A/R$ has trivial radical as $R$ is the radical of $A$. By Lemma~\ref{lem=insolvable}, $A/R$ has a normal arc-transitive subgroup $I/R$ such that $B/R\leq  I/R$ and $(B/R,I/R)\cong (G,T)=(M_{11},M_{12})$ or $(A_{n-1},A_n)$ with $n\geq 8$ and $n\di 2^4\cdot 3^6$.

Note that $I\unlhd A$. Since $A$ is $2$-arc-transitive, $I$ is arc-transitive. Let $C=C_I(R)$. Then $C\unlhd I$ and $C\cap R\leq Z(C)$. Since $B=R\times G\leq I$, we have $G\leq C$, and since $C/C\cap R\cong CR/R\unlhd I/R\cong T$, we have $C\cap R=Z(C)$, $C/Z(C)\cong T$ and $I=CR$. Furthermore, $C'/C'\cap Z(C)\cong C'Z(C)/Z(C)=(C/Z(C))'=C/Z(C)\cong  T$. Thus, $Z(C')=C'\cap Z(C)$, $C=C'Z(C)$ and $C'/Z(C')\cong T$. It follows that $C'=(C'Z(C))'=C''$, and so $C'$ is a covering group of $T$. Since $C/C'$ is abelian, $G\leq C'$.

Recall that $(G,T)=(M_{11},M_{12})$ or $(A_{n-1},A_n)$ with $n\geq 8$ and $n\di 2^4\cdot 3^6$. By \cite[Theorem 5.1.4]{Kleidman} and \cite[pp.31]{Atlas}, the Schur multiplier $\Mult(A_n)=\mz_2$ for $n\geq 8$ and $\Mult(M_{12})=\mz_2$. Then $T$ has a unique double cover, denoted by $2.T$.

Suppose $Z(C')=\mz_2$. Then $C'$ is the unique double cover of $T$, that is, $C'=2.T$. Let  $(G,T)=(M_{11},M_{12})$. Then $|C'|=2|M_{12}|$, and since $\G$ is $G$-regular, $|C'_v|=24$  for $v\in V(\G)$. In particular, $3\di |I_v|$, and $I$ is $2$-arc-transitive by the arc-transitivity of $I$. Since $C\unlhd I$ and $C'$ is characteristic in $C$, we have $C'\unlhd I$ and so $C'$ is arc-transitive. Similarly, $C'$ is $2$-arc-transitive as $3\di |C'_v|$, and by Proposition~\ref{prop=stabilizer}, $C'_v\cong S_4$. The quotient graph $\G_{Z(C')}$ is a connected tetravalent $(C'/Z(C'),2)$-arc-transitive graph with stabilizer isomorphic to $C'_v$. Since $C'_v\cong S_4$ and $C'/Z(C')\cong M_{12}$, this is impossible by Lemma~\ref{A12}.
Now let $(G,T)=(A_{n-1},A_n)$ with $n\geq 8$.  Then $C'=2.A_n$. Since $Z(C')\unlhd I$ and $I/R\cong T$, we have $Z(C')\leq R$, and since $B=G\times R$, we have $G\times Z(C')\leq C'$. Then $G\times Z(C')$ is a subgroup of index $n$ of $C'$ isomorphic to $A_{n-1}\times \mz_2$.  This is impossible by Proposition~\ref{prop=covering group}.

Thus, $Z(C')=1$. It follows that $C'\cong T$ and $G\leq  C'\unlhd I$. Since $|I|=|I/R||R|=|T||R|=|C'||R|$ and $C'\cap R=1$, we have $I=C'\times R$ because $C'\unlhd I$, and so $C'$ is characteristic in $I$. Since $I\unlhd A$ and $A$ is $2$-arc-transitive, $A$ has a normal arc-transitive non-abelian simple subgroup $C'\cong T$ containing $G$ and $(G,T)=(M_{11},M_{12})$ or $(A_{n-1},A_n)$ with $n\geq 8$ and $n\di 2^4\cdot 3^6$. Since $\G$ is $G$-regular, $|C'_v|=|C'_v|/|G_v|=|C'|/|G|=|T:G|$, and by Proposition \ref{prop=stabilizer}, $(G,T)=(M_{11},M_{12})$ or $(A_{n-1},A_n)$ with $n=2^2\cdot 3$, $2^3\cdot 3$, $2^2\cdot3^2$, $2^3\cdot3^2$, $2^4\cdot3^2$, $2^4\cdot3^3$ or $2^4\cdot3^6$. Furthermore, $(G,T)\not=(A_{11},A_{12})$ by Lemma~\ref{A12}. It follows that $(G,T)$ is listed in Table~\ref{table=1}.
\hfill\qed

Now, we are ready to prove Theorem~\ref{theo=main}.

\medskip

\f {\bf The proof of Theorem \ref{theo=main}:} Let $G$ be a non-abelian simple group and $\G$ a connected tetravalent $2$-arc-transitive $G$-regular graph with   $v\in V(\G)$. Then $G_v=1$ and $|G|=|V(\G)|$. Set $A=\Aut(\G)$ and let $R$ be the radical of $A$. By Lemma \ref{lem=insolvable}, the theorem is true for $R=1$. Thus, we may assume $R\not=1$.

Set $B=RG$. Then $G\cap R=1$ and so $|B|=|R||G|$. Since $\G$ is $G$-regular,  $B=GB_v$ and $|B|=|G||B_v|$. It follows that $|R|=|B_v|$, and by Proposition~\ref{prop=stabilizer}, $|R|\di 2^4\cdot 3^6$.

Suppose that $R$ has one or two orbits on $V(\G)$. Since $\G$ is a connected tetravalent $G$-regular graph, $|R|=|R_v||v^R|=|R_v||G|$ or $|R_v||G|/2$. Since $|R|\di 2^4 \cdot 3^6$, the non-abelian simple group $G$ is a $\{2,3\}$-group,  which is impossible.

Thus, $R$ has at least three orbits. If $B=R\times G$ then the theorem is true by Lemma \ref{lem=T normal}. Now assume $B\not=R\times G$, and to finish the proof, we aim to derive contradictions.

Since $|R|\di 2^4\cdot 3^6$, we may write $|R|=2^m\cdot3^k$, where $0\leq m\leq 4$ and
$0\leq k\leq 6$. Since $R$ is soluble, there exists a series of principle subgroups of $B$:
\begin{center}
$B>R=R_{s}>\cdots >R_1>R_0=1,$
\end{center}
such that $R_i\unlhd B$ and $R_{i+1}/R_i$ is an elementary abelian $r$-group with $0\leq i\leq s-1$ and $r=2$ or $3$. Let $|R_{i+1}/R_i|=r^{\ell_i}$. Then $\ell_i\leq m\leq 4$ if $r=2$, and $\ell_i\leq k\leq 6$ if $r=3$. Note that $G\leq B$ has a natural action on $R_{i+1}/R_i$ by conjugation.

Since $B\neq R\times G$, there exists $0\leq j\leq s-1$ such that $GR_j=G\times R_j$ and $GR_{j+1}\not=G\times R_{j+1}$. If $G$ acts trivially on $R_{j+1}/R_j$, then $[GR_j/R_j,R_{j+1}/R_j]=1$. Since $GR_j/R_j\cong G$ is simple, $(GR_j/R_j)\cap (R_{j+1}/R_j)=1$, and since $|GR_{j+1}/R_j|=|GR_{j+1}/R_{j+1}||R_{j+1}/R_j|=|G||R_{j+1}/R_j|=|GR_j/R_j||R_{j+1}/R_j|$, we have $GR_{j+1}/R_j=GR_j/R_j\times R_{j+1}/R_j$. In particular, $GR_j\unlhd GR_{j+1}$ and so $G\unlhd GR_{j+1}$ because $GR_j=G\times R_j$ implies that $G$ is characteristic in $GR_j$. It follows that $GR_{j+1}=G\times R_{j+1}$, a contradiction. Thus, $G$ acts non-trivially on $R_{j+1}/R_j$, and the simplicity of $G$ implies that the action is faithful. Since $R_{j+1}/R_j$ is an elementary abelian group of order $r^{\ell_j}$, we have $G\leq \GL(\ell_j,r)$, where $\ell_j\leq m\leq 4$ if $r=2$ and $\ell_j\leq k\leq 6$ if $r=3$.

By Proposition \ref{prop=atlesst3orbits}, $R$ is semiregular on $V(\G)$, and the quotient graph $\G_{R}$ is a connected tetravalent $(A/R,2)$-arc-transitive graph with $A/R\leq \Aut(\G_{R})$. Moreover, $\G_R$ is $B/R$-vertex-transitive, and $|B/R|=|V(\G_R)||(B/R)_\a|$ for $\a\in V(\G_R)$. Since $\G$ is $G$-regular, $|V(\G_R)|=|V(\G)|/|R|=|G|/|R|$, and since $B/R\cong G$, we have $|G|=|B/R|=|V(\G_R)||(B/R)_\a|=|G|/|R|\cdot |(B/R)_\a|$. It follows that $|R|\di |G|$ and $|(B/R)_\a|=|R|$.

Since $R$ is the largest normal soluble subgroup of $A$, $A/R$ has trivial radical, and since $B/R\cong G$, Lemma~\ref{lem=insolvable} implies that either $B/R\unlhd A/R$, or $A/R$ has a normal arc-transitive subgroup $I/R$ such that $B/R\leq I/R$ and $(B/R,I/R) \cong (G,T)=(M_{11},M_{12})$ or $(A_{n-1},A_n)$ with $n\geq 8$ and $n\di 2^4\cdot 3^6$.

\medskip
\f {\bf Case 1:} $B/R\unlhd A/R$.

In this case, $B\unlhd A$ and $B_v\unlhd A_v$. Since $A$ is $2$-arc-transitive, $B$ is arc-transitive and  $\G_R$ is $B/R$-arc-transitive. Thus, $4\di |B_v|$ and   $\G_R\cong \Cos(B/R,(B/R)_\a, (B/R)_\a g (B/R)_\a)$ for some feasible $g$.
Recall that $|R|=|B_v|=|(B/R)_\a|=2^m\cdot 3^k$, $|R|\di |G|$ and $G\leq \GL(\ell_j,r)$, where $\ell_j\leq m\leq 4$ if $r=2$ and $\ell_j\leq k\leq 6$ if $r=3$.

By Proposition \ref{prop=stabilizer}, $A_v\cong A_4, \mz_3\times A_4, S_4,  \mz_3\rtimes S_4, S_3\times S_4, \AGL(2,3)$ or $\mz_3^5 \rtimes \GL(2,3)$.

Suppose $A_v\cong A_4$, $\mz_3\times A_4$, $S_4$, $\mz_3\rtimes S_4$ or $ S_3\times S_4$. Since $|R|=|B_v|$, we have $|R|=2^m\cdot 3^k$, where  $m\leq 4$ and $k\leq 2$.  Since $\GL(2,2)$ and $\GL(2,3)$ are soluble and $G\leq \GL(\ell_j,r)$, we have either $m=4$ and $G\leq \GL(4,2)=\PSL(4,2)$, or $m=3$ and $G=\PSL(3,2)$.

If $m=4$ and $G\leq \PSL(4,2)$  then $2^4\di |B_v|$ and $A_v=S_3\times S_4$. Since $B_v\unlhd A_v$, we have $B_v=S_3\times S_4$, and by Proposition~\ref{prop=stabilizer}, $|(B/R)_\a|=|B_v|=2^4\cdot 3^2$ implies that  $(B/R)_\a=S_3\times S_4$. This is impossible because $B/R\cong G\leq \PSL(4,2)$ and $\PSL(4,2)$ has no subgroup isomorphic to $S_3\times S_4$ by MAGMA.
If $m=3$ and $G=\PSL(3,2)$ then $A_v=S_4$, $\mz_3\rtimes S_4$ or $ S_3\times S_4$. Since $B_v\unlhd A_v$ and $2^3\di |B_v|$, we have $|(B/R)_\a|=|B_v|=2^3\cdot 3, 2^3\cdot 3^2, 2^4\cdot 3$ or $2^4\cdot 3^2$, and by Proposition~\ref{prop=stabilizer}, $(B/R)_\a\cong S_4, \mz_3\rtimes S_4$ or $S_3\times S_4$. By Atlas \cite[pp.3]{Atlas}, $(B/R)_\a\cong S_4$ as $B/R\cong \PSL(3,2)$. By MAGMA, $(B/R)_\a$ has two conjugacy classes in $B/R$, and for both classes, there are no feasible $g$, a contradiction.

Suppose $A_v\cong\rm AGL(2,3)$. Since $4\di |B_v|$ and $B_v\lhd A_v$, by MAGMA $B_v=A_v$, and so $|R|=2^4\cdot 3^3$. Since $3^3\nmid |\GL(4,2)|$, we have $G\leq \GL(3,3)$, and by MAGMA, $G=\PSL(3,3)$.  Since $|(B/R)_\a|=|R|=2^4\cdot 3^3$, by Proposition~\ref{prop=stabilizer} $(B/R)_\a\cong \rm AGL(2,3)$. By MAGMA, there are two conjugacy classes isomorphic to $\AGL(2,3)$ in $\PSL(3,3)$, and  for each conjugacy class, there is no feasible $g$, a contradiction.

Suppose $A_v\cong \mz_3^5 \rtimes \rm GL(2,3)$. Since $4\di |B_v|$ and $B_v\unlhd A_v$, by the Remark of Proposition~\ref{prop=stabilizer}, $B_v=A_v$ and $|R|=|(B/R)_\a|=2^4\cdot 3^6$. Thus, $B$ is $7$-arc-transitive and $G\cong B/R$ has a subgroup of order $2^4\cdot 3^6$. Since $2^4\cdot3^6\nmid |\GL(4,2)|$, we have  $G\nleq \GL(4,2)$. Thus, $G\leq \GL(\ell_j,3)$ with $\ell_j\leq k\leq 6$. By MAGMA, for $1\leq k\leq 4$, $\GL(k,3)$  have no simple subgroup $G$ with a subgroup of order $2^4\cdot 3^6$. It follows that $\ell_j=5$ or $6$ with $r=3$, and hence $|R_{j+1}/R_j|=3^5$ or $3^6$, where $0\leq j\leq s-1$.

Let $j\not=s-1$. Then $|R/R_{s-1}|=3$ or $2^t$ for $1\leq t\leq 4$, and so $|R_{s-1}|=2^4\cdot 3^5$ or $2^{4-t}\cdot 3^6$. Since $G$ is simple and $G\nleq \GL(4,2)$, $G$ acts trivially on $R/R_{s-1}$ by conjugation, implying $[GR_{s-1}/R_{s-1},R/R_{s-1}]=1$. It follows that $GR/R_{s-1}=GR_{s-1}/R_{s-1}\times R/R_{s-1}$ and hence $GR_{s-1}\unlhd GR=B$. Since $B$ is $7$-arc-transitive, $GR_{s-1}$ is arc-transitive, and since $G$ is regular, $GR_{s-1}$ has a stabilizer of order $2^4\cdot 3^5$ or $2^{4-t}\cdot3^6$, contrary to Proposition~\ref{prop=stabilizer}.

Let $j=s-1$. Then $|R/R_{s-1}|=3^5$ or $3^6$, and $|R_{s-1}|=2^4\cdot 3$ or $2^4$. Furthermore, $GR_{s-1}=G\times R_{s-1}$ and $B=GR\not=G\times R$. Clearly, $R=R_{s-1}P$ for a Sylow $3$-subgroup $P$ of $R$. Let $C=C_B(R_{s-1})$. Then $C\unlhd B$ and $G\leq C$. If $G=C$ then $G\lhd B$ and so $B= G\times R$, a contradiction. Thus, $G$ is a proper subgroup of $C$, so that $C_v\not=1$. Since $B$ is $7$-arc-transitive, $C$ is arc-transitive and by the Remark of Proposition \ref{prop=stabilizer}, $C_v=B_v$, that is, $B=C$.
In particular, $[P,R_{s-1}]=1$. Since $R=PR_{s-1}$, $P$ is normal in $R$ and so characteristic. This implies $P\unlhd B$ and so $GP\leq B=C$. It follows $[GP,R_{s-1}]=1$, and since $B=GR=(GP)R_{s-1}$, we have $GP \unlhd B$. Since $B$ is $7$-arc-transitive, $GP$ is arc-transitive, forcing $4\di |(GP)_v|$, and since $G$ is regular, $|(GP)_v|=|P|$ and thus $4\di |P|$, contrary to the fact that $P$ is a Sylow $3$-subgroup of $R$.

\medskip
\f {\bf Case 2:} $A/R$ has a normal arc-transitive subgroup $I/R$ such that $B/R\leq I/R$ and $(B/R,I/R) \cong (G,T)= (M_{11},M_{12})$ or $(A_{n-1},A_n)$ with $n\geq 8$ and $n\di 2^4\cdot 3^6$.

In this case, $I\unlhd A$ and $I$ is arc-transitive. Since $|T|=|I/R|=|V(\G_R)||(I/R)_\a|=|G|/|R|\cdot |(I/R)_\a|$,  we have $|(I/R)_\a|=|R||T|/|G|$. By Proposition~\ref{prop=stabilizer}, $|R||T|/|G|$ is a divisor of $2^4\cdot 3^6$. Recall that $|R|=|B_v|=|(B/R)_\a|=2^m\cdot 3^k$, $|R|\di |G|$ and $G\leq \GL(\ell_j,r)$, where $\ell_j\leq m\leq 4$ if $r=2$ or $\ell_j\leq k\leq 6$ if $r=3$.

Suppose $(G,T)=(M_{11},M_{12})$. Then $|T|/|G|=2^2\cdot 3$. Since $|R||T|/|G|\di 2^4\cdot3^6$, we have $|R|\di 2^2\cdot3^5$, and then $|R|\di |G|$ implies that $|R|\di 2^2\cdot3^3$. Since $\GL(2,2)$ is soluble, we have  $G\leq \GL(3,3)$, which is impossible because $11\di |M_{11}|$ and $11\nmid |\GL(3,3)|$.

Suppose $(G,T)=(A_{n-1},A_n)$ with $n\geq 8$ and $n\di 2^4\cdot 3^6$.  If $n\geq 12$ then $5^2\di |G|$, which is impossible because $G\leq \GL(4,2)$ or $\GL(6,3)$ but $5^2\nmid |\GL(6,3)|$ and $5^2\nmid \GL(4,2)$. It follows \textcolor[rgb]{1,0,0}{that} $8\leq n<12$, and since $n\di 2^4\cdot 3^6$, we have  $(G,T)=(A_7,A_8)$ or  $(A_8,A_9)$.

Suppose $(G,T)=(A_7,A_8)$. Then $|T|/|G|=2^3$. Since $|R||T|/|G|\di 2^4\cdot3^6$, we have $|R|\di 2\cdot3^6$, and $|R|\di |G|$ implies that $|R|\di 2\cdot3^2$. Thus, $G\leq \GL(1,2)$ or $\GL(2,3)$, yielding that $G$ is soluble, a contradiction.

Suppose $(G,T)=(A_8,A_9)$. Then $|T|/|G|=3^2$. Since $|R||T|/|G|\di 2^4\cdot3^6$, we have $|R|\di 2^4\cdot 3^4$ and $|R|\di |G|$ implies $|R|\di 2^4\cdot3^2$, that is, $|R|=2^m\cdot 3^k$ with $m\leq 4$ and $k\leq 2$.
Since $\GL(2,3)$ is soluble, $G\leq \GL(4,2)$, and hence $G=\GL(4,2)$ and $m=4$ as $G=A_8\cong \GL(4,2)$. Since
$|(I/R)_\a|=|R||T|/|G|$, we have $2^4\cdot 3^2\di|(I/R)_\a|$, and since $|I/R|=|A_9|=2^6\cdot3^4\cdot5\cdot7$, Proposition \ref{prop=stabilizer} implies that $(I/R)_\a\cong S_3\times S_4$ or $\rm AGL(2,3)$. By the arc-transitivity of $I$, $\G_R$ is $I/R$-arc-transitive, and so $\G_R=\Cos(I/R,(I/R)_\a , (I/R)_\a t (I/R)_\a)$ for some feasible $t$.
By MAGMA, $I/R\cong A_9$ has no subgroup isomorphic to $\rm AGL(2,3)$ and $I/R$ has two conjugacy classes isomorphic to $S_3\times S_4$. By taking a given $(I/R)_\a$ in each conjugacy class, computation shows that there is no feasible $t$, a contradiction.

\hfill\qed

To end the paper, we give examples to show that the pair $(G,T)=(A_{23},A_{24})$ in Theorem~\ref{theo=main} can happen.

\begin{exam}\label{NonNormalExample} Let $G=A_{23}$ and $T=A_{24}$. Define $x,y,z,w,g\in T$ as following:

\begin{enumerate}
\item[]$x=(1, 2)(3, 7)(4, 10)(5, 13)(6, 15)(8, 12)(9, 19)(11, 18)(14, 22)(16, 20)(17,
    24)(21, 23)$,\\
$y=(1, 3)(2, 7)(4, 8)(5, 9)(6, 18)(10, 12)(11, 15)(13, 19)(14, 20)(16, 22)(17,
    23)(21, 24)$,\\
$z=(1, 4, 6)(2, 8, 11)(3, 12, 15)(5, 17, 16)(7, 10, 18)(9, 21, 20)(13, 23, 14)(19,
    24, 22)$,\\
$w=(1, 5)(2, 9)(3, 13)(4, 16)(6, 17)(7, 19)(8, 20)(10, 22)(11, 21)(12, 14)(15,
    23)(18, 24)$, \\
$g_1=(1, 5)(2, 10)(3, 14)(4, 17)(6, 16)(7, 11)(8, 18)(9, 22)(12, 13)(15, 23)(19,
    21)(20, 24)$,\\
$g_2=(1, 5)(2, 10)(3, 8)(4, 16)(6, 15)(7, 19)(9, 22)(11, 12)(13, 20)(14, 18)(17,
    24)(21, 23)$, \\
$g_3=(1, 5)(2, 9)(3, 13)(4, 16)(6, 15)(7, 14)(8, 20)(10, 19)(11, 17)(12, 23)(18,
    22)(21, 24)$, \\
$g_4=(1, 2)(3, 12)(4, 8)(5, 9)(6, 10)(7, 19)(11, 23)(13, 22)(14, 18)(15, 20)(16,
    21)(17, 24)$. \\
\end{enumerate}

\rm {By MAGAMA\cite{magma}, $H=\langle x,y,z,w\rangle\cong S_4$, $T=\langle H,g_i\rangle$, $|H:H\cap H^{g_i}|=4$ $(1\leq i\leq 4)$ and $H$ is regular on $\{1,2,\cdots, 24\}$. Thus, $T=GH$ with $G\cap H=1$. Define $\G_i=\Cos(T,H,Hg_iH)$ with $1\leq i\leq 4$. Then $\G_i$ is a connected tetravalent $(T,2)$-arc-transitive $G$-regular graph, where $G$ and $T$ are viewed as groups of automorphisms of $\G_i$ by right multiplication. By Theorem~\ref{theo=main}, $T\unlhd \Aut(\G_i)$ with $1\leq i\leq 4$. Again by MAGMA, $\Aut(T,H,HgH)\cong \tilde{H}$, where $\tilde{H}$ is the automorphism group of $T$ induced by conjugate of elements in $H$. Thus, $\Aut(\G_i)=T$ ($1\leq i\leq 4$) by~\cite[Lemma~2.10]{WFZ}.}
\end{exam}

\medskip
\f {\bf Acknowledgement:} This work was supported by the National Natural Science Foundation of China (11571035, 11231008) and by the 111 Project of China (B16002).

\newpage
\begin{center}{\bf \Large Appendix}\end{center}

\vskip 0.2cm

Let $G$ be a finite group and $H\leq G$. Let $\G=\Cos(G,H,HgH)$ be a connected tetravalent $G$-arc-transitive graph. Then $g$ can be chosen as a $2$-element such that
$g^2\in H$, $\langle H,g\rangle=G$, $|H:H\cap H^g|=4$, and such a $g$ is called feasible. There are many places in the paper to compute feasible $g$ for given $G$ and $H$ and to compute their full automorphism groups of the corresponding coset graphs. Here we provide a computer program based on MAGMA language by taking $(G,H)=(M_{12}, A_4)$ as an example:

\medskip
\f load``M12";
   B:=Subgroups(G);

\medskip
\f PG:=[ ]; //possible graphs

\medskip
\f for i in [1..\#B] do if IsIsomorphic(\rm Alt(4),B[i]'subgroup) then
H:=B[i]'subgroup; print i;

\medskip
\f D:=[ ]; // feasible elements

\medskip
\f for $g$ in $G$ do
   if IsDivisibleBy $(2^6,\rm Order(g))$ and Order(sub $\langle G| H,g^2\rangle$) eq Order$(H)$ and \#$(H*g*H)$ eq Order$(H)$*4 and Order(sub $\langle G| H, g\rangle$) eq Order$(G)$ then Include ($\sim\rm D,g$);
  end if;
 end for;

\medskip
\f \#D;

\medskip
\f if \#D ne 0 then for j in [1..\#D] do
c:=D[j]; HcH:=\{\};

\medskip
\f for $t$ in $H$ do
for $h$ in $H$ do
Include $(\sim \rm HcH, t*c*h)$;
end for; end for;

\medskip

\f Vj:=\{\}; Ej:=\{\};

\f for $t$ in $G$ do for $s$ in HcH do
T1:=\{\}; T2:=\{\}; for $h$ in $H$ do
Include $(\sim \rm T1, h*t)$; Include $(\sim \rm T2, h*s*t)$; end for; Include $(\sim \rm Vj, T1)$; Include $(\sim \rm Ej, \{T1,T2\})$; end for; end for;

\medskip
\f PGj:=Graph $\langle \rm Vj | \rm Ej \rangle$;
Include($\sim \rm PG,PGj$);

\f end for; end if; end if; end for;

\medskip
\f NPG:=[]; // non-isomorphic possible graphs

\medskip
\f NPG:=[PG[1]];

\f for k in [1..\#PG] do
p:=0;
      for m in [1..\#NPG] do
        if IsIsomorphic(PG[k], NPG[m]) then
         p:=p+0;
        else p:=p+1;
     end if;
     end for;

\f if p eq \#NPG then
  Include($\sim \rm NPG,PG[k]$);
end if;
end for;

\medskip

\f Graph:=[ ]; //$M_{12}$-arc-transitive $M_{11}$-regular graph

\medskip
\f for t in [1..\#NPG] do  A:=AutomorphismGroup $(\rm NPG[t])$; S:=Subgroups(A);

\f for n in [1..\#S] do
     if IsTransitive(S[n]'subgroup) and Order(S[n]'subgroup) eq 7920
and IsSimple(S[n]'subgroup) then Include($\sim\rm Graph,NPG[t]$);

\f  print ``We find a connected tetravalent $M_{12}$-arc-transitive  $M_{11}$-regular graph";

\f print ``The automorphism group A:";  print A;

\f print ``The radical of A:"; print Radical(A);

\f end if; end for; end for;

\medskip
\f \#Graph;

\end{document}